\documentclass[11pt]{amsart}
\usepackage{amssymb}
\usepackage{times}
\usepackage{graphicx}

\theoremstyle{plain}

\numberwithin{equation}{section}

\newcommand{\calD}{\mathcal{D}}

\newcommand{\calH}{\mathcal{H}}
\newcommand{\calI}{\mathcal{I}}

\newcommand{\calL}{\mathcal{L}}

\newcommand{\calS}{\mathcal{S}}

\newcommand{\bbF}{\mathbb{F}}
\newcommand{\bbC}{\mathbb{C}}

\newcommand{\bbP}{\mathbb{P}}
\newcommand{\bbQ}{\mathbb{Q}}
\newcommand{\bbR}{\mathbb{R}}

\newcommand{\bbZ}{\mathbb{Z}}

\newcommand{\la}{\langle}
\newcommand{\ra}{\rangle}

\def\GL{{\text{GL}}}
\def\SL{{\text{SL}}}

\def\Pic{{\text{Pic}}}

\def\O{{\text{O}}}

\def\rank{{\text{rank}}}

\def\mod{{\text{mod}}}
\def\dim{{\text{dim}}}
\def\trace{{\text{trace}}}

\def\Sp{{\text{Sp}}}
\def\SO{{\text{SO}}}

\def\I{{\text{I}}}
\def\II{{\text{II}}}

\def\IV{{\text{IV}}}

\begin{document}
\title [Igusa quartic and Borcherds products] {Igusa quartic and Borcherds products}
\author{Shigeyuki Kond\=o}
\address{Graduate School of Mathematics, Nagoya University, Nagoya,
464-8602, Japan}
\email{kondo@math.nagoya-u.ac.jp}
\thanks{Research of the author is partially supported by
Grant-in-Aid for Scientific Research S-22224001, Japan}

\begin{abstract}
By applying Borcherds' theory of automorphic forms on bounded 
symmetric domains of type $\IV$, 
we give a 5-dimensional linear system of automorphic forms of weight $6$
on Igusa quartic $3$-fold which induces an $\mathfrak{S}_6$-equivariant rational map of degree $16$ from Igusa quartic to Segre cubic.  In particular we have a rational self-map of
Igusa quartic of degree $16$.
\end{abstract}
\maketitle

\section{Introduction}

The purpose of this paper is to give an application of the theory of automorphic forms on bounded symmetric domains 
of type $\IV$ due to Borcherds \cite{B1}, \cite{B2}.  We consider Igusa quartic 3-fold $\calI$ given by
\begin{equation}\label{igusa}
\sum_i x_i = \Bigl(\sum_i x_i^2 \Bigr)^2 - 4\Bigl( \sum_i x_i^4 \Bigr) = 0 \subset \bbP^5,
\end{equation}
where $(x_1:\cdots : x_6)$ is a homogenous coordinate of $\bbP^5$.
 It is classically known (Baker \cite{Ba}, Chap.V, \ Dolgachev \cite{D}) that Igusa quartic is the dual variety of
Segre cubic 3-fold $\calS$ defined by
\begin{equation}\label{segre}
\sum_i x_i =  \sum_i x_i^3 = 0.
\end{equation}
The symmetry group $\mathfrak{S}_6$ of degree 6 naturally acts on $\calI$ and $\calS$ 
as automorphisms.   
Igusa quartic $\calI$ is isomorphic to the Satake compactification 
$\overline{\mathfrak{H}_2/\Gamma_2(2)}$
of the quotient of the Siegel 
upper half plane $\mathfrak{H}_2$ of degree two by the $2$-congrunce subgroup 
$\Gamma_2(2)$ of $\Gamma_2=\Sp(4,\bbZ)$  (Igusa \cite{I}, page 397; 
also see van der Geer \cite{vG}).   The natural action of $\mathfrak{S}_6 (\cong \Gamma_2/\Gamma_2(2)$) on $\overline{\mathfrak{H}_2/\Gamma_2(2)}$ coincides
with the above one on $\calI$.
On the other hand, let 
$M=U(2)^{\oplus 2}\oplus A_1(2)$ be the transcendental lattice of a generic Kummer surface associated with a smooth
curve of genus 2 and let $\calD(M)$ be a bounded symmetric domain of type $\IV$ and of
dimension $3$ associated with $M$.  Then 
it follows from Gritsenko, Hulek \cite{Gr}, \cite{GHu} that 
$\overline{\mathfrak{H}_2/\Gamma_2(2)}$ is isomorphic to the Baily-Borel compactification $\overline{\calD(M)/\Gamma_M}$ where $\Gamma_M$ 
is a subgroup of the orthogonal group $\O(M)$ such that $\O(M)/\Gamma_M \cong \mathfrak{S}_6$.  This isomorphism is $\mathfrak{S}_6$-equivariant.

In this paper, 
by applying Borcherds' theory,
we present holomorphic automorphic forms $\Psi_{10}, \Psi_{30}, \Psi_{24}$ on $\calD(M)$ of weight $10, 30, 24$ which coincide with the Siegel modular forms with Humbert surfaces $H_1, H_4, H_5$ as zero divisors, respectively (Theorem \ref{Geer}). 
Moreover we give a 5-dimensional linear system of holomorphic automorphic forms on 
$\calD(M)$ of weight $6$ 
which induces the linear system $\calL$ of cubics on $\calI$ given by
\begin{equation}\label{15functions}
(x_i - x_j)(x_k-x_l)(x_m-x_n)\quad  (\{i,j,k,l,m,n\} = \{1,\ldots, 6\})
\end{equation}
(Theorem \ref{restr2}).
The linear system $\calL$ gives
an $\mathfrak{S}_6$-equivariant rational map of degree $16$ from Igusa 
quartic to Segre cubic (Theorem \ref{Igor}).  Thus we have a rational self-map of Igusa quartic of degree 
$16$ (Corollary \ref{selfmap}).  Note that Mukai \cite{M} recently found a holomorphic self-map of Igusa quartic of degree $8$.
On the other hand, Segre cubic is an arithmetic quotient of a 3-dimensional complex ball (see Hunt \cite{H}).  This complex ball can be naturally embedded into a bounded symmetric domain of type $\IV$ and of dimension $6$.
Recently, by applying Borcherds theory, the author \cite{K3} gives a 5-dimensional
space of automorphic forms on the complex ball which defines the dual map from Segre cubic to Igusa quartic.

We use an idea of Allcock, Freitag \cite{AF} in which they gave an embedding of the moduli space of marked cubic surfaces into $\bbP^9$ by applying Borcherds' theory.  For a given 
lattice $L$ of signature $(2,n)$, we consider vector-valued modular forms with respect to
the Weil representation of $\SL(2,\bbZ)$ on the group ring $\bbC[L^*/L]$ where $L^*$ is the dual of $L$.  There are
two types of liftings of vector-valued modular forms both of which give automorphic forms
on the bounded symmetric domain $\calD(L)$ of type $\IV$ associated with $L$.  One is called {\it Borcherds product} or {\it multiplicative lifting} which is an automorphic form with known zeros and poles.  Another one is called {\it additive lifting} which is an automorphic form with respect to the subgroup of the orthogonal group $\O(L)$ acting trivially on $L^*/L$. 
Borcherds gives explicit formulae for the Fourier coefficients of the additive lifting in terms
of the Fourier coefficients of the vector-valued modular form.

We apply Borcherds' theory to the lattice 
$N=U(2)^{\oplus 2} \oplus A_1^{\oplus 2}$ of signature $(2,4)$ instead of $M$ 
because $M$ has {\it odd} rank $5$ which makes a difficulty of calculations.  
We can embed $M$ into $N$ as a primitive sublattice which induces an embedding
of the domain $\calD(M)$ into $\calD(N)$.  
We remark that $\calD(N)$ is the period domain of $K3$ surfaces associated with six lines on $\bbP^2$ (see Matsumoto, Sasaki, Yoshida \cite{MSY}).
We construct automorphic forms $\Phi_4, \Phi_{10}, \Phi_{30}, \Phi_{48}$ on $\calD(N)$ of weight $4, 10, 30, 48$ as Borcherds products (Theorem \ref{borcherds}).
By restricting them to $\calD(M)$, we get Siegel modular forms 
$\Psi_{10}, \Psi_{30}, \Psi_{24}$ mentioned as above ($\Phi_4$ vanishes on $\calD(M)$).  For example, 
$\Psi_{10}$ is the product $\prod \theta_m^2(\tau)$ of the square of even theta constants.
We remark that the {\it Borcherds product} $\Phi_4$ is also obtained by {\it additive lifting} (Remark \ref{1dim}).   Such example was already given by the author in the case of the
moduli space of Enriques surfaces (\cite{K1}, Remark 4.7).
On the other hand, to each 2-dimensional isotorpic subspace in $A_N=N^*/N$,
we associate an automorphic form of weight 10 as additive lifting.  By restricting them,
we get fifteen automorphic forms on $\calD(M)$ of weight 6 corresponding to
fifteen functions given in (\ref{15functions}) (Theorem \ref{restr2}).

The plan of this paper is as follows.  
In \S \ref{sec2}, we recall a theory of lattices.
Section \ref{six-lines} is devoted to the theory of periods of $K3$ surfaces
which are double covers of $\bbP^2$ branched along six lines, due to Matsumoto, Sasaki,
Yoshida \cite{MSY}.  In section \ref{AQ}, we recall a description of the Igusa quartic as an 
arithmetic quotient of a bounded symmetric domain of type $\IV$.  Moreover we study the boundary components and Heegner divisors (= Humbert surfaces) on the Igusa quartic.  
In section \ref{weil}, we recall the Weil representation of $\SL(2,\bbZ)$ on the group ring
$\bbC [N^*/N]$ and calculate its character. We study the 5-dimensional subspace appeared in the Weil representation.  The result will be used to construct additive liftings in \S \ref{GBL}. 
In section \ref{BP}, by using Borcherds products, we show that 
there exist holomorphic automorphic forms on $\calD(N)$ of
weight $4, 10, 30, 48$ with known zeros, and in section \ref{GBL}, by using additive liftings, we give a 5-dimensional space of automorphic forms
on $\calD(N)$.  Finally, in \S \ref{G}, we discuss automorphic forms on Igusa quartic.

\medskip
\emph{Acknowledgments}:
The author thanks Klaus Hulek and Matthias Sch\"utt for valuable conversations, and Igor Dolgachev for 
discussions in Schiermonnikoog 2014.  In particular the proof of Theorem \ref{Igor} is due to Dolgachev.

\medskip

\section{Preliminaries}\label{sec2}

A {\it lattice}  is  a free abelian group $L$ of  finite rank equipped with 
a non-degenerate symmetric integral bilinear form $\langle, \rangle : L \times L \to \bbZ$. 
 For $x \in L\otimes \bbQ$, we call $x^2 =x\cdot x$ the {\it norm} of $x$.
For a lattice $L$ and a rational number
$m$, we denote by $L(m)$ the free $\bbZ$-module $L$ with the $\bbQ$-valued bilinear form obtained from the bilinear form of $L$ by multiplication with $m$. The signature of a lattice is the signature of the real quadratic space $L\otimes \bbR$. 
A lattice is called {\it even} if $\langle x, x\rangle \in 2\bbZ$ for all $x\in L$. 
  
We denote by $U$ the even unimodular lattice of signature $(1,1)$, 
and by $A_m, \ D_n$ or $\ E_k$ the even {\it negative} definite lattice defined by
the Cartan matrix of type $A_m, \ D_n$ or $\ E_k$ respectively.  
For an integer $m$, we denote by $\langle m\rangle$ the lattice of rank 1 generated by a vector with norm $m$.  
We denote by $L\oplus M$ the orthogonal direct sum of lattices $L$ and $M$, and by $L^{\oplus m}$ the orthogonal direct sum of $m$-copies of $L$.

Let $L$ be an even lattice and let $L^* ={\rm Hom}(L,\bbZ)$.  
We denote by $A_L$ the quotient $L^*/L$ which is called the {\it discriminant group} of $L$, and define maps
$$q_L : A_L \to \bbQ/2\bbZ, \quad b_L : A_L \times A_L \to \bbQ/\bbZ$$
by $q_L(x+L) = \langle x, x\rangle\ {\rm mod}\ 2\bbZ$ and
$b_L(x+L, y+L) = \langle x, y\rangle \ {\rm mod} \ \bbZ$.  
We call $q_L$ the {\it discriminant quadratic form} of $L$ and $b_L$ the {\it discriminant bilinear form}.  
A lattice is called {\it 2-elementary} if its discriminant group is a 2-elementary abelian group.  
We denote by $u$, $v$ or $q_{\pm}$ the discriminant quadratic form of 2-elementary 
lattice $U(2)$, $D_4$ or $\langle \pm 2\rangle$ respectively.  For any 2-elementary 
lattice $L$, the discriminant form $q_L$ is a direct sum of $u, v, q_{\pm}$.
An even 2-elementary lattice $L$ is called type $\I$ if $q_L$ is a direct sum
of $u$ and $v$, and type $\II$ if otherwise.  It is known that 
the isomorphism class of an even indefinite 2-elementary lattice is determined by
its signature, the rank of $A_L$ and its type $\I$ or $\II$.


Let ${\rm O}(L)$ be the orthogonal group of $L$, that is, the group of isomorphisms of $L$ preserving the bilinear form.
Similarly ${\rm O}(q_L)$ denotes the group of isomorphisms of $A_L$ preserving $q_L$.
There is a natural map
\begin{equation}\label{nat}
{\rm O}(L) \to {\rm O}(q_L)
\end{equation}
whose kernel is denoted by $\widetilde{\O}(L)$.
For more details we refer the reader to Nikulin \cite{N1}.

\section{Six lines on $\bbP^2$ and $K3$ surfaces}\label{six-lines}

Let $\ell_1,\ldots, \ell_6$ be ordered six lines on $\bbP^2$ 
in general position, that is, no three lines meet at one point.
Let $\bar{X}$ be the
double cover of $\bbP^2$ branched along the sextic $\ell_1 +\cdots + \ell_6$.
Then $\bar{X}$ has $15$ ordinary nodes over the point
$p_{ij} = \ell_i \cap \ell_j$.  Let $X$ be the minimal resolution of $\bar{X}$ which is 
a $K3$ surface.  Obviously $X$ contains $15$ mutually disjoint smooth rational curves
$E_{ij}$ which are the exceptional curves over $p_{ij}$, and six smooth
rational curves $\tilde{\ell}_i$ $(1\leq i\leq 6)$ which are the proper transforms of $\ell_i$.
Denote by $\Pic(X)$ the Picard lattice of $X$ and by $S_X$ the smallest primitive sublattice of $\Pic(X)$ containing $21$ smooth rational curves $E_{ij}, \ \tilde{\ell}_i$.

\subsection{Proposition}(Matsumoto, Sasaki, Yoshida \cite{MSY})\label{pic}
{\it The lattice 
$S_X$ is a primitive sublattice of $H^2(X,\bbZ)$ of signature $(1,15)$, 
$S_X^*/S_X \cong (\bbZ/2\bbZ)^6$ and $q_{S_X} \cong u \oplus u\oplus q_{+}
\oplus q_{+}$.  The group $\O(q_{S_X})$ is isomorphic to ${\mathfrak S}_6 \times \bbZ/2\bbZ$, where ${\mathfrak S}_6$ is the symmetric group  of degree $6$.
The natural map $\O(S_X) \to \O(q_{S_X})$ is surjective.
}
\begin{proof}  The assertion follows from Corollary 2.1.6 and Proposition 2.8.2 in \cite{MSY}.
Here we give an another proof by using Nikulin's lattice theory.  First note that $S_X$ is the invariant sublattice of $H^2(X,\bbZ)$ under the
action of the covering transformation of $X \to \bbP^2$.
It follows from Nikulin \cite{N2}, Theorem 4.2.2 that 
$S_X$ is an even 2-elementary lattice of signature
$(1,15)$ and with $q_{S_X} = u\oplus u\oplus q_+ \oplus q_+$.
Note that there exists a subgroup 
$F \cong (\bbZ/2\bbZ)^5$ of $A_{S_X}$
on which the restriction of $q_{S_X}$ has values in $\bbZ/2\bbZ$, that is, $q_{S_X}|F$
is a quadratic form of dimension $5$ over $\bbF_2$.   Moreover
$q_{S_X}|F$ contains a radical $\langle \kappa \rangle \cong \bbZ/2\bbZ$, that is,
$\kappa$ is perpendicular to all elements in $F$ with respect to $b_{S_X}$,  
and  $b_{S_X}$ 
induces a symplectic form on $F/\langle \kappa \rangle$ of dimension 4 over $\bbF_2$.  Thus the orthogonal group $\O (q_{S_X})$ is isomorphic to 
$\Sp(4,\bbF_2) \times \bbZ/2\bbZ \cong {\mathfrak S}_6 \times \bbZ/2\bbZ$
where $\bbZ/2\bbZ$ is generated by the involution changing two components $q_+\oplus q_+$ and acting trivially on $u\oplus u$.
The surjectivity of the natural map
$\O(S_X) \to \O(q_{S_X})$ follows from \cite{N1}, Theorem 3.6.3.
\end{proof}

We call $X$ {\it generic} if $S_X = \Pic(X)$.
Let $T_X$ be the orthogonal complement of $\Pic(X)$ in $H^2(X,\bbZ)$ which is called
the {\it transcendental lattice} of $X$.  Also we denote by $N_X$ the orthogonal complement of $S_X$.  It is known that $q_{S_X} = - q_{N_X} $
(e.g. Nikulin \cite{N1}, Corollary 1.6.2) and hence $q_{N_X} \cong u\oplus u \oplus q_-\oplus q_-$.
Since $N_X$ is a $2$-elementary lattice of signature $(2,4)$, the isomorphism class of $N_X$ is determined by its signature and $q_{N_X}$
(Nikulin \cite{N1}, Theorem 3.6.2).  Thus we have
$$N_X \cong U(2)^{\oplus 2} \oplus A_1^{\oplus 2} \cong A_1(-1)^{\oplus 2}\oplus A_1^{\oplus 4}.$$
We denote by $N$ an abstract lattice of signature $(2,4)$ and with 
$q_N=u\oplus u \oplus q_-\oplus q_-$.  If $X$ is generic, then 
$T_X \cong N$.  We denote by $\kappa_N \in A_N$ the radical corresponding to $\kappa$
(see the proof of Proposition \ref{pic}).

An elementary calculation shows the following Lemma.

\subsection{Lemma}\label{types}
{\it The discriminant group $A_{N}$ consists of the following $64$ vectors}:

\smallskip
Type $(00) : \alpha = 0, \ \#\alpha =1$;
\smallskip

Type $(0) :  \alpha \not= 0,\ q_{N}(\alpha) = 0,\ \# \alpha = 15;$
\smallskip

Type $(1) : \alpha \not= \kappa_N,\ q_{N}(\alpha) = 1,\ \# \alpha = 15;$
\smallskip

Type $(10) : \alpha= \kappa_N,\  \# \alpha = 1;$ 
\smallskip

Type $(1/2) : q_{N}(\alpha) = 1/2,\ \# \alpha = 12;$
\smallskip

Type $(3/2) : q_N(\alpha)=3/2,\  \# \alpha = 20.$ 
\smallskip

\medskip
Define
$$\calD(N) = \{ [\omega] \in \bbP(N\otimes \bbC) \ : \ \la \omega, \omega\ra = 0, \ \la \omega, \bar{\omega}\ra > 0\}$$
where $\bar{\omega}$ is the complex conjugate of $\omega$.  
It is known that $\calD(N)$ is a disjoint union of two copies of
a bounded symmetric domain of type $\IV$ and of dimension $4$.  
We denote by $\Gamma_N$ the group $\widetilde{\O}(N)$ which acts properly discontinuously
on $\calD(L)$.  It is known that the 
Baily-Borel compactification of the quotient $\calD(N)/\Gamma_N$ is the coarse moduli space of ordered six lines on $\bbP^2$
(\cite{MSY}).

Now we define the Heegner divisors on $\calD(N)$.  
Fix a vector  $\alpha \in A_N$ with $q_N(\alpha)\not= 0$
and a negative rational number $n$.  
For $r\in N^*$ with $r^2<0$, we denote by $r^{\perp}$ the hyperplane in $\calD(N)$ perpendicular to $r$.
We define a {\it Heegner divisor} $\calH(N)_{\alpha,n}$ by
$$\calH(N)_{\alpha, n} = \sum_r r^{\perp}$$
where $r$ runs through all vectors $r$ in $N^*$ satisfying 
$r \ \mod \ N = \alpha$ and $\la r, r\ra = n$.
For simplicity, we denote by $\calH(N)_{\alpha}$ 
the Heegner divisor $\calH(N)_{\alpha, n}$ for $n = -1, -1/2$ or $-3/2$ according to $q_N(\alpha) = 1, 3/2$ or $1/2$, respectively.  
Also we denote by $\calH(N)_{1}, \calH(N)_{3/2}$ or $\calH(N)_{1/2}$
the union of all $\calH(N)_{\alpha}$ where $\alpha$ runs through all vectors $\alpha$ 
with $q_N(\alpha) = 1$ ($\alpha \not=\kappa_N$), $3/2$ or $1/2$ respectively.  
The geometric meaning of these Heegner divisors is known.  For example, 
A generic point in $\calH(N)_{1}$ corresponds to six lines such that three points
$p_{ij}, p_{kl}, p_{mn}$ are collinear where $\{i,j,k,l,m,n\} =\{1,\ldots,6\}$.
For more details we refer the reader to \cite{LPS}, Theorem 3.6.

\subsection{Reflections}\label{reflection}
Let $r$ be a $(-4)$-vector in $N$ with $r/2 \in N^*$.
Then the reflection $s_r$ defined by
\begin{equation}\label{ref2}
s_r(x) = x - {2\la x, r\ra \over \la r, r \ra} r = x + \la x, r/2\ra r \quad (x \in N)
\end{equation}
is contained in $\O(N)$.  The reflection $s_r$ induces a reflection
$t_{\alpha}$ on $A_N$ associated with $\alpha = r/2 \ \mod \ N$ defined by
\begin{equation}\label{ref}
t_{\alpha}(\beta) = \beta + 2b_N(\beta, \alpha)\alpha \quad (\beta \in A_N).
\end{equation}

\section{Igusa quartic and a bounded symmetric domain of type $\IV$ }\label{AQ}


Let ${\mathfrak H}_2$ be the Siegel upper half plane of degree two and let 
$\Gamma_2(2)$ be the principal 
2-congurence subgroup of $\Gamma_2 = \Sp(4, \bbZ)$.  We denote by $\overline{\mathfrak H_2/\Gamma_2(2)}$ the Satake compactification of
the quotient ${\mathfrak H}_2/\Gamma_2(2)$.  Igusa \cite{I} showed that $\overline{\mathfrak H_2/\Gamma_2(2)}$ can be embedded into $\bbP^4$ by using
theta constants, whose image is a quartic hypersurface $\calI$ given
by the equation (\ref{igusa}) called {\it Igusa quartic} (Igusa gave a different form.  See \cite{vG}).  
The boundary of the compactification consists of fifteen 1-dimensional components 
and fifteen 0-dimensinal components which correspond to fifteen lines and fifteen points on the Igusa quartic.
By definition ${\mathfrak H}_2/\Gamma_2(2)$ is the moduli space of principally polarized abelian surfaces with a level 2-structure.
 
On the other hand it is known that ${\mathfrak H}_2$ is isomorphic to a bounded symmetric domain of type $\IV$ and of dimension 3
as bounded symmetric domains (e.g. \cite{vG}).  Put
\begin{equation}\label{kummer-trans}
M = U(2)^{\oplus 2} \oplus A_1(2).
\end{equation}
Then $M$ is an even lattice of signature $(2,4)$ and is isomorphic to the transcendental lattice of a generic Kummer surface associated with a smooth curve of genus two.
Define
$$\calD(M) = \{ [\omega] \in \bbP(M\otimes \bbC) \ : \ \la \omega, \omega\ra = 0, \ \la \omega, \bar{\omega}\ra > 0\}$$
which is a disjoint union of two copies of
a bounded symmetric domain of type $\IV$ and of dimension 3.  
The quotient $\calD(M)/\O(M)$ is birational to the moduli space of Kummer surfaces associated with a smooth curve of genus 2.

The discriminant group $A_M =M^*/M$ is isomorphic to 
$(\bbZ/2\bbZ)^4 \oplus\bbZ/4\bbZ$.
Let $A^{(2)}_M$ be the 2-elementary subgroup of $A_M$.
Then $A^{(2)}_M  \cong (\bbZ/2\bbZ)^5$ and the restriction of the discriminant form $q_M$ on $A^{(2)}_M$  
has
a radical $\bbZ/2\bbZ$.  We denote by $\kappa_M$ the generator of the radical.
If $a$ is a generator of the component $A_1(2)$ of the decomposition (\ref{kummer-trans})  of $M$, then
$\kappa_M = a/2\ \mod \ M$.  The discriminant bilinear form $b_M$ induces a symplectic form on $A^{(2)}_M /\langle \kappa_M\rangle$ of dimension $4$ over $\bbF_2$.  Thus
the orthogonal group $\O (q_M)$ is isomorphic to 
$\Sp(4,\bbF_2) \times \bbZ/2\bbZ \cong {\mathfrak S}_6 \times \bbZ/2\bbZ$
where $\bbZ/2\bbZ$ is generated by the involution $-1_{A_M}$.
Put
$$\widetilde{\O}(M) = \{ \gamma \in \O(M) \ : \ \gamma \mid A_M = 1\}, \quad  \Gamma_M = \{ \gamma \in \O(M) \ : \ \gamma \mid A_M = \pm 1\}.$$
Obviously we have two exact sequences
$$1 \longrightarrow \Gamma_M \longrightarrow \O(M)
\longrightarrow {\mathfrak S}_6  \longrightarrow 1, \quad 1 \longrightarrow \widetilde{\O}(M) \longrightarrow \Gamma_M
\longrightarrow \bbZ/2\bbZ \longrightarrow 1.$$  

Gritsenko, Hulek \cite{GHu}, \S 1 gave an explicit correspondence between  $\Sp_4(\bbZ)$ and $\SO(3,2)_{\bbZ}$.
First we remark that in the paper \cite{GHu} they considered the lattice $M(-1/2) = U^{\oplus 2}\oplus A_1(-1)$ instead of $M$.  
However $M$ is obtained from 
$M(-1/2)$ by multiplying the bilinear form by $-2$, and hence $\O(M) \cong \O(M(-1/2))$. 
They gave an isomorphism 
$$\Psi : \Gamma_2 \to \SO(M) \cap \O^+(M)$$
explicitly, where $\SO(M)$ is the special orthogonal group and $\O^+(M)$ is the subgroup of $\O(M)$
preserving a component of $\calD(M)$.
Since $M$ has odd rank,  $-1_M$ and $\SO(M)$ generate $\O(M)$.
Since $-1_M$ acts trivially on 
$\calD(M)$ and $\O(M)$ interchanges two components of $\calD(M)$, we have an isomorphism
${\mathfrak H}_2/\Gamma_2 \cong \calD(M)/\O(M).$
Note that $M^*/M \cong M/2M$.  By using this fact and the explicit  isomorphism given in \cite{GHu}, we  see that
the image of the principal 2-congurence subgroup $\Gamma_2(2)$ is 
contained in $\widetilde{\O}(M) \cap \SO(M)$.
Since $-1_M \in \Gamma_M$ represents $-1_{A_M}$ and acts trivially on $\calD(M)$, we have an isomorphism
$${\mathfrak H}_2/\Gamma_2(2) \cong \calD(M)/\Gamma_M.$$
Now we conclude that Igusa quartic $\calI$ is isomorphic to the Baily-Borel compactification $\overline{\calD(M)/\Gamma_M}$
of the quotient $\calD(M)/\Gamma_M$.


Next we study boundaries and Heegner divisors on $\overline{\calD(M)/\Gamma_M}$.

\subsection{Lemma}\label{types2}
{\it The discriminant group $A_{M}/\{\pm 1_{A_M}\}$ consists of the following $48$ vectors}:

Type $(00)$ : $\alpha = 0, \ \#\alpha =1$;
\smallskip

Type $(0)$ : $ \alpha \not= 0,\ q_{M}(\alpha) = 0,\ \# \alpha = 15,$
\smallskip

Type $(1)$ : $q_{M}(\alpha) = 1,\ \alpha \not= \kappa_M,\  \# \alpha = 15,$
\smallskip

Type $(10)$ : $\alpha= \kappa_M,\  \# \alpha = 1,$ 
\smallskip

Type $(3/4)$ : $q_{M}(\alpha) = 3/4,\ \# \alpha = 6,$
\smallskip

Type $(7/4)$ : $q_M(\alpha)=7/4,\  \# \alpha = 10.$ 
\begin{proof}
The assertion follows from a direct calculation.
We remark that the involution $-1_{A_M}$ acts trivially on the 2-elementary subgroup $A^{(2)}_M $ of $A_M$.
In particular $-1_{A_M}$ fixes all vectors of type $(00), (0), (1), (10)$, and
$$-1_{A_M}(\alpha) = \alpha + \kappa_M$$
if $\alpha$ is of type $(3/4), (7/4)$.
\end{proof}



\subsection{Boundary components}\label{boundary20}
We call a subgroup $T$ of $M$ an {\it isotoropic sublattice} if the symmetric bilinear form
vanishes on $T$.  Since $M$ has the signature $(2,4)$, the rank of an isotoropic sublattice
is at most $2$.  Similary we define an {\it isotoropic subspace} of $A_M$ as a subgroup
on which the discriminant quadratic form $q_M$ vanishes.  The dimension of an isotoropic
subspace is also at most $2$.

It is known that $\overline{\calD(M)/\Gamma_M} (\cong \overline{{\mathfrak H}_2/\Gamma_2(2)})$ has fifteen $0$-dimensional boundary components and
fifteen $1$-dimensional boundary components.
$0$-dimensional (resp. 1-dimensional) boundary components bijectively correspond to primitive isotoropic sublattices of rank $1$ (resp. of rank 2) 
in $M$ modulo $\Gamma_M$.  
 A primitive isotoropic sublattice of rank 1 (resp.  isotoropic
sublattice of rank 2) in $M$ determines a non-zero isotoropic vector (resp. 1-dimensional
isotoropic subspace) in $A_M$.  For example, if $\langle e_1, e_2\rangle$ is a primitive  isotoropic sublattice of $M$ generated by $e_1, e_2$, then 
$\langle e_1/2 \ \mod \ M, \ e_2/2\ \mod \ M\rangle$ is an isotoropic subspace in $A_M$.

\subsection{Lemma}\label{boundary}
{\it The $0$-dimensional $($resp. $1$-dimensional$)$ boundary components correspond to non-zero isotoropic vectors $($resp. $1$-dimensional isotoropic subspaces$)$ in $A_M$.}
\begin{proof}
Since $\Gamma_M$ acts trivially on isotoropic vectors in $A_M$, it suffices to see that
there exists exactly $15$ non-zero isotoropic vectors and 15 isotoropic subspaces in $A_M$.
The first assertion follows from Lemma \ref{types2}.  Moreover we see that for each non-zero
isotoropic vector $\alpha \in A_M$, there are 7 non-zero isotropic vectors (including $\alpha$) perpendicular to $\alpha$.  This implies that there are three isotoropic subspaces containing $\alpha$.  Since there are 15 non-zero isotoropic vectors and each 1-dimensional
isotoropic subspace contains three non-zero isotoropic vectors, the number of
1-dimensional isotoropic subspaces is $(15 \times 3)/3 = 15$.
\end{proof}

\subsection{Remark}\label{boundary2}
The incidence relation between 15 $0$-dimensional boundary components and 
15 $1$-dimensional boundary components is called {\it $(15)_3$-configuration} because
each $0$-dimensional boundary component is contained in exactly three $1$-dimensional boundary components
and each $1$-dimensional boundary component contains exactly three $0$-dimensional
boundary componets (e.g. see \cite{vG}).

\subsection{Heegner divisors}\label{heegner2}
Next
we define Heegner divisors on  $\calD(M)$ as those on $\calD(L)$.  
Let $r\in M^*$ with $r^2 <0$.
Denote by $r^{\perp}$ the hyperplane in $\calD(M)$ orthogonal to $r$.
Fix a vector  $\alpha \in A_M$ with $q_M(\alpha)\not= 0$,
$\alpha\not= \kappa_M$ and a negative rational number $n$.  
We define a {\it Heegner divisor} $\calH(M)_{\alpha,n}$ by
$$\calH(M)_{\alpha, n} = \sum_r r^{\perp}$$
where $r$ runs through all vectors $r$ in $M^*$ satisfying $r \ \mod \ M = \alpha$ and $\la r, r\ra = n$.
For simplicity we denote by $\calH(M)_{\alpha}$ the Heegner divisor $\calH(M)_{\alpha, n}$ for  $n = -1, -5/4$ or $-1/4$ according to $q_M(\alpha) = 1, 3/4$ or $7/4$ respectively.
We also denote by $\calH(M)_{1}, \calH(M)_{3/4}$ or $\calH(M)_{7/4}$
the union of all $\calH(M)_{\alpha}$ where $\alpha$ runs through all vectors of type $(1), (3/4)$ or $(7/4)$ respectively.
The image of a Heegner divisor in $\calD(M)/\Gamma_M$ is also called a {\it Heegner divisor}.

In Gritsenko, Hulek \cite{GHu}, Lemma 3.2, they proved that any two vectors in $M(1/2)^*$ with the same norm and the same image in $A_{M(1/2)}$ are conjugate under the action of $\O(M(1/2))$.  
It follows  that all $(-1)$-vectors $r \in M^*$ with $r \ \mod \ M$ being of type  $(1)$
are conjugate under the action of $\O(M)$.  The same statement holds for $(-5/4)$- or $(-1/4)$-vectors $r \in M^*$ with $r \ \mod \ M$ being of type  $(3/4)$ or
$(7/4)$ respectively. 
Therefore 
$\calH(M)_{1}/\Gamma_M, \calH(M)_{3/4}/\Gamma_M$ or 
$\calH(M)_{7/4}/\Gamma_M$ has 
exactly fifteen, six or ten irreducible components $\calH(M)_{\alpha}/\Gamma_M$ where
$\alpha \in A_M/\{\pm 1\}$ is of type $(1), (3/4)$ or $(7/4)$ respectively (see Lemma \ref{types2}).  

\subsection{Remark}\label{Humbert}

In the theory of moduli of abelian surfaces, Heegner divisors are called {\it Humbert surfaces} (e.g. \cite{vG}, \cite{GHu}).
Let us compare Heegner divisors and Humbert surfaces.  Recall that $M=U(2)^{\oplus 2} \oplus A_1(2)$ is obtained from
$M(1/2) =U^{\oplus 2} \oplus A_1$ by multiplying the bilinear form by 2.
In the notation as in \cite{GHu}, the {\it Humbert surface $H_{\Delta}$ of the discriminant }
$\Delta$ is the image of the Heegner divisor
$\calH_{-\Delta/2}$ on $\calD(M(1/2))$ because they consider the lattice of signature 
$(3,2)$ and hence we should take the opposite 
sign (see the definition on page 476 in \cite{GHu}).  The Heegner divisor $\calH_{-\Delta/2}$ on $\calD(M(1/2))$ corresponds to
the Heegner divisor $\calH_{-\Delta/4}$ on $\calD(M)$.  Thus the closure of 
$\calH(M)_{7/4}/\Gamma_M, \calH(M)_{1}/\Gamma_M$ or 
$\calH(M)_{3/4}/\Gamma_M$ in Baily-Borel compactification of 
$\calD(M)/\Gamma_M$ is equal to
the Humbert surface $H_1$, $H_4$ or $H_5$ given in \cite{GHu}, \cite{vG}, respectively.

\section{The Weil representation}\label{weil}

In this section we recall the Weil representation associated with the lattice 
$N = U(2) \oplus U(2)\oplus A_1\oplus A_1$ given in
\S \ref{six-lines}, and calculate its character. 
The following Table \ref{table1} means that for each vector $u \in A_N$ of given type, 
$m_j$ is the number of vectors $v \in A_N$ 
of given type with $\langle u, v \rangle = j/2$.

\begin{table}[h]
\[
\begin{array}{rlllllllllllllllllllllll}
u& 00&00&00&00&00&00&0&0&0&0&0&0& 1&1&1&1&1&1\\
v&00&0&1&10&3/2&1/2&00&0&1&10&3/2&1/2&00&0&1&10&3/2&1/2\\
m_0&1&15&15&1&20&12&1&7&7&1&12&4&1&7&7&1&8&8\\
m_1&0&0&0&0&0&0&0&8&8&0&8&8&0&8&8&0&12&4\\
u&10&10&10&10&10&10&3/2&3/2&3/2&3/2&3/2&3/2&1/2&1/2&1/2&1/2&1/2&1/2\\
v&00&0&1&10&3/2&1/2&00&0&1&10&3/2&1/2&00&0&1&10&3/2&1/2\\
m_0&1&15&15&1&0&0&1&9&6&0&10&6&1&5&10&0&10&6\\
m_1&0&0&0&0&20&12&0&6&9&1&10&6&0&10&5&1&10&6\\
\end{array}
\]
\caption{}
\label{table1}
\end{table}


Let 
$
T =
\begin{pmatrix}1&1
\\0&1
\end{pmatrix},\ 
S =
\begin{pmatrix}0&-1
\\1&0
\end{pmatrix}.
$
Then $S$ and $T$ generate $\SL(2,\bbZ)$.
Let $\rho$ be the {\it Weil representation} of $\SL(2, \bbZ)$ on the group ring $\bbC[A_N]$ defined by

\begin{equation}\label{weil2}
\rho(T)(e_{\alpha}) = e^{\pi \sqrt{-1} q_N(\alpha)} e_{\alpha}, \quad
\rho(S)(e_{\alpha}) = {\sqrt{-1}\over 8} \sum_{\beta \in A_N} 
e^{-2\pi \sqrt{-1} b_N( \beta, \alpha)} e_{\beta}.
\end{equation}


By definition and Table \ref{table1}, we see that $\rho(S^2)(e_{\alpha}) = - e_{\alpha}.$
The action of $\SL(2,\bbZ)$ on $\bbC[A_N]$ factorizes to the one of $\SL(2,\bbZ/4\bbZ)$.  
The conjugacy classes of $\SL(2, \bbZ/4\bbZ)$ consist of
$\pm E, \pm S, \pm T, \pm T^{2}, ST, (ST)^2$.  
Let $\chi_{i}$ $(1 \leq i \leq 10)$ be the characters of irreducible
representations of $\SL(2, \bbZ/4\bbZ)$.  
One can easily compute the character table of 
$\SL(2, \bbZ/4\bbZ)$ by using GAP \cite{GAP08}.
For the convenience of the reader we give the
character table (Table \ref{table2}) of $\SL(2, \bbZ/4\bbZ)$.

\begin{table}[h]
\[
\begin{array}{rlllllllllllllllllllllll}
 & E&-E&S&-S&T&-T&T^2&-T^2&ST&(ST)^2\\
\chi_1&1&1&1&1&1&1&1&1&1&1\\
\chi_2&1&1&-1&-1&-1&-1&1&1&1&1\\
\chi_3&1&-1&\sqrt{-1}&-\sqrt{-1}&\sqrt{-1}&-\sqrt{-1}&-1&1&-1&1\\
\chi_4&1&-1&-\sqrt{-1}&\sqrt{-1}&-\sqrt{-1}&\sqrt{-1}&-1&1&-1&1\\
\chi_5&2&2&0&0&0&0&2&2&-1&-1\\
\chi_6&2&-2&0&0&0&0&-2&2&1&-1\\
\chi_7&3&3&1&1&-1&-1&-1&-1&0&0\\
\chi_8&3&3&-1&-1&1&1&-1&-1&0&0\\
\chi_9&3&-3&-\sqrt{-1}&\sqrt{-1}&\sqrt{-1}&-\sqrt{-1}&1&-1&0&0\\
\chi_{10}&3&-3&\sqrt{-1}&-\sqrt{-1}&-\sqrt{-1}&\sqrt{-1}&1&-1&0&0\\
\end{array}
\]
\caption{}
\label{table2}
\end{table}

\subsection{Lemma}\label{char2}
{\it Let $\chi$ be the character of the Weil representation of
$SL(2, \bbZ/4\bbZ)$ on $\bbC[A_N]$.
Let 
$$\chi = \sum_{i=1}^{10} m_{i} \chi_{i}$$ 
be
the decomposition of $\chi$ into irreducible characters.  Then
$$\chi = \chi_3+ 5\chi_4+5\chi_6 + 6\chi_9+10\chi_{10}.$$}
\begin{proof}
By definition (\ref{weil2}) and Table \ref{table1}, 
we see that $\trace (E) = 2^6$, $\trace (-E) = -2^6$, $\trace (S) = 0$, $\trace(-S)=0$,
$\trace(T)= -8\sqrt{-1}$, $\trace(-T)=8\sqrt{-1}$, $\trace(T^2)=0$, $\trace(-T^2)=0$, $\trace(ST) =-1$ and $\trace((ST)^2)=1$.  The assertion now follows from the Table \ref{table2}.
\end{proof}
\smallskip

\subsection{Definition}\label{5-dim}
Let $W$ (resp. $W_0$) be the subspace in $\bbC [A_N]$ on which the character of $\SL(2,\bbZ)$ is given by $5\chi_4$ (resp. $\chi_3$).  Note that the action of $\O(q_N)$ on $\bbC[A_N]$ 
commutes with the action of $\SL(2,\bbZ)$.  Therefore $\O(q_N)$ acts on $W$ and $W_0$.
In the section \ref{GBL} we will construct a 5-dimensional space of automorphic forms on $\calD(N)$ associated with $W$ (For $W_0$, see Remark \ref{1dim}).

\subsection{Definition}\label{invariant}
Let $I$ be a 2-dimensional isotropic
subspace  of $A_N$ with respect to $q_N$.  Note that $I$ is a maximal isotoropic subspace. Let $V$ be the subspace of $A_N$ generated by $I$ and $\kappa_N$.  
Take a vector $\alpha_0 \in A_N$ satisfying 
$q_N(\alpha_0) = 3/2$ and $b_N(\alpha_0, c) = 0$ for any $c \in I$.  Note that
$\alpha_0$ is unique modulo $V$ because $I^{\perp}/V = \bbF_2$.
Define
$$M_{+} = \{ \alpha_0 + c : c \in I\}, \quad M_{-} = \{ \alpha_0 + c + \kappa_N : c \in I\},$$
and
$$\theta_V = \sum_{\beta \in M_{+}} e_{\beta} - \sum_{\beta \in M_{-}}
e_{\beta} \in \bbC[A_N].$$
This definition is the same as the one given in the case of the moduli space of plane quartic curves
in \cite{K2}.

\subsection{Lemma}\label{key}
\

(i)  \ $\rho(S)(\theta_{V}) = -\sqrt{-1} \theta_{V}$ and $\rho(T)(\theta_{V}) = -\sqrt{-1} \theta_{V}.$
{\it In particular $\theta_V$ is contained in $W$.}
\smallskip

(ii) \ {\it  For $\beta \in V$ with $q_N(\beta) = 1$, 
$t_{\beta}(\theta_{V}) = -\theta_{V}$ where $t_{\beta}$ is the reflection associated 
with $\beta$.}

\begin{proof}
(i)  If $\beta \in M_{\pm}$, then $q_N(\beta)=3/2$, and hence 
$\rho(T)(\theta_{V}) = -\sqrt{-1} \theta_{V}.$  Next by definition (\ref{weil2}),
\begin{equation}\label{kon}
\rho(S)(\theta_{V}) = {\sqrt{-1} \over 8} \sum_{\delta} \left(\sum_{\beta\in M_+}
e^{-2\pi\sqrt{-1}\ b_N(\delta, \beta)} - \sum_{\beta\in M_-} e^{-2\pi\sqrt{-1}\ 
b_N(\delta, \beta)}\right) e_{\delta}.
\end{equation}
We denote by ${\sqrt{-1} \over 8}\cdot c_{\delta}$ the coefficient of $e_{\delta}$ in the equation (\ref{kon}).  
If $\delta \in M_+$, then
$b_N(\delta, \beta) = 1/2$ for $\beta \in M_+$ and 
$b_N(\delta, \beta) \in \bbZ/2\bbZ$ for $\beta \in M_-$.  
Hence $c_{\delta} = -2^2-2^2=-2^3$.
Similary if $\delta \in M_-$, then $c_{\delta} = 2^3$.

Now assume $\delta \notin M_{\pm}$.  If $\delta \in V$, we can easily see that $c_{\delta} =0$.
Hence we assume $\delta \notin V$.
First consider the case 
$b_N(\delta, \kappa_N) \in \bbZ/2\bbZ$.
Since $V^{\perp}=V$, there exists $\gamma \in V$ such that 
$b_N(\gamma, \delta) \notin \bbZ/2\bbZ$.  In this case 
$I = \delta^{\perp} \cap I \cup \{ \gamma + a : a \in \delta^{\perp} \cap I \}$.  
This implies that
$$\sum_{\beta\in M_+}
e^{-2\pi\sqrt{-1} \ b_N(\delta, \beta)} = \sum_{\beta\in M_-} e^{- 2\pi\sqrt{-1}\ b_N(\delta, \beta)}=0.$$
Finally if $b_N(\delta, \kappa_N) \notin \bbZ/2\bbZ$, then 
$\delta = \alpha_0 + \delta'$ and 
$b_N(\delta', \kappa_N) \in \bbZ/2\bbZ$.  Then this case reduces to the previous case.

(ii) Let $\beta \in V$.  Then $\beta = c + \kappa_N$, $c \in I$.  
If $c' \in I$, then $\langle \beta, \alpha_0 + c' \rangle = 1/2$.  Therefore
the reflection $t_{\beta}$ defined by the equation (\ref{ref}) interchanges $M_+$ and $M_-$ and hence the assertion follows.
\end{proof}

\subsection{Lemma}\label{totally}
{\it There exist exactly fifteen $2$-dimensional isotropic subspaces in $A_N$.}
\begin{proof}
Recall that each non-zero isotropic vector $\alpha \in A_N$, 
there exist exactly 7 non-zero isotropic vectors 
(including $\alpha$) perpendicular to $\alpha$ (see Table \ref{table1}).  It follows that there are three maximal totally isotropic subspaces containing
$\alpha$.   Since the number of non-zero isotropic vectors is 15, the number of maximal totally isotropic subspace is
$15 \times 3/3 =15$.
\end{proof}

Thus we have 15 vectors $\theta_V$ in $W$.

\subsection{Lemma}\label{irred} 
{\it As a $\O(q_N) (\cong {\mathfrak S}_6 \times {\bbZ}/2{\bbZ})$ module, $W$ is irreducible.}
\begin{proof}
It is well known that there are no irreducible representations of $\mathfrak{S}_6$ of degree
$2,3,4$.
If $W'$ is an irreducible representation of ${\mathfrak S}_6$ and $\dim \ W' \geq 2$, then $\dim \ W' \geq 5$.  Hence it suffices to see that there are 
no 1-dimensional invariant subspaces under the action of ${\mathfrak S}_6$.  
Assume that $W$ is a direct sum of 1-dimensional representations.
A direct calculation shows that there exist five linearly independent vectors $\theta_{V_1},\ldots, \theta_{V_5}$ where $V_i$ $(i=1,\ldots, 5)$ is a subspace of $A_N$ generated by a 2-dimensional isotropic subspace and $\kappa_N$.   It follows from Lemma \ref{key}, (ii) that $W$ is a direct sum of
alternating representations.  In particular any 1-dimensional subspace is invariant under the action of ${\mathfrak S}_6$. 
However any  vector $\theta_V$ as above is not invariant under the action of ${\mathfrak S}_6$.  This is a contradiction.
\end{proof}

In the section \ref{GBL} we will construct a 5-dimensional space of automorphic forms on $\calD(N)$ associated with $W$.

\subsection{Correction}\label{}
The proof of Lemma 9 in \cite{K3} is not complete.  We should add a sentence "There exist five lineraly independent vectors $\nu_{\alpha_0^1}, \ldots, \nu_{\alpha_0^5}$."  Then the same proof as in Lemma \ref{irred} holds
in this case, too.


\section{Borcherds products}\label{BP}

{\it Borcherds products} are automorphic forms on $\calD(N)$ 
whose zeros and poles lie on Heegner divisors.  First of all we recall the definition of automorphic forms.  Define
$$\widetilde{\calD}(N) =\{ \omega \in N\otimes {\bf C} \ : \ \la \omega, \omega\ra = 0, \ \la \omega, \bar{\omega} \ra > 0\}.$$
Then the canonical map $\widetilde{\calD}(N) \to {\calD}(N)$ is a ${\bf C}^*$-bundle.
Let $\Gamma$ be a subgroup of $\O(N)$ of finite index.
A holomorphic function
$$\Phi : \widetilde{\calD}(N) \to {\bf C}$$
is called a holomorphic automorphic form of weight $k$ with respect to $\Gamma$ on $\calD(N)$ if $\Phi$ is homogeneous of degree $-k$, that is,
$\Phi (c\cdot \omega) = c^{-k}\Phi (\omega)$ for $c\in {\bf C}^*$, and is invariant under $\Gamma$.

Let $\gamma$ be a representation of  $\SL(2,\bbZ)$ on $\bbC[A_N]$. 
A {\it vector-valued modular form of weight $k$ and of type $\gamma$} is a holomorphic map 
$$f = \{ f_{\alpha}\}_{\alpha\in A_N} : H^+ \to \bbC [A_N]$$
satisfying 
$$f(A\tau) = (c\tau + d)^{k} \gamma(A) f(\tau)$$
where $A =
\begin{pmatrix}a& b
\\c &d
\end{pmatrix} \in \SL(2,\bbZ)$.  We assume that $f$ is holomorphic at cusps.

In this section, we will show that there exists a holomorphic automorphic form
of weight $4$, $10$, $30$ or $48$ whose zero divisor is the Heegner divisor 
$\calH(N)_{\kappa_N}$, 
$\calH(N)_{3/2}$, $\calH(N)_{1}$ or $\calH(N)_{1/2}$ respectively.
To show the existence of such 
Borcherds products, we introduce the {\it obstraction space}
consisting of all vector-valued modular forms $\{ f_{\alpha} \}_{\alpha \in A_N}$ of weight $3=(\rank(N)/2)$ with respect to
the dual representation $\rho^*$ of the Weil representation $\rho$ given in (\ref{weil2}):

\begin{equation}\label{dualweil2}
\rho^*(T)(e_{\alpha}) = e^{-\pi \sqrt{-1}\langle \alpha, \alpha
\rangle} e_{\alpha}, \quad
\rho^*(S)(e_{\alpha}) = -{\sqrt{-1}\over 8} \sum_{\beta \in A_N} 
e^{2\pi \sqrt{-1} \langle \beta, \alpha \rangle} e_{\beta}.
\end{equation}
In other words,
\begin{equation}\label{dual}
f_{\alpha}(\tau + 1) = e^{-\pi \sqrt{-1}\ \la \alpha, \alpha \ra}f_{\alpha}(\tau), \quad
f_{\alpha}(-1/\tau) = -{\sqrt{-1} \tau^{3} \over 8} 
\sum_{\beta} e^{2\pi \sqrt{-1} \ \la \alpha, \beta \ra} f_{\beta}(\tau).
\end{equation}

\noindent
We will apply the next theorem to show the existence of such Borcherds products.

\subsection{Theorem}\label{freitag}(Borcherds \cite{B2}, Freitag \cite{F}, Theorem 5.2)
{\it A linear combination
$$\sum_{\alpha \in A_N, \ n<0} c_{\alpha, n} \calH(N)_{\alpha, n} \quad (c_{\alpha, n} \in \bbZ)$$
of Heegner divisors
is the divisor of an automorphic form on $\calD(N)$ 
of weight $k$ if for every cusp form 
$$f = \{f_{\alpha}(\tau)\}_{\alpha\in A_N}, 
\quad f_{\alpha}(\tau) = \sum_{n \in \bbQ} a_{\alpha, n} e^{2\pi \sqrt{-1} n \tau}$$ 
in the obstruction space, the relation
$$\sum_{\alpha \in A_N, \ n<0} a_{\alpha, -n/2} c_{\alpha, n} = 0$$
holds.  In this case the weight $k$ is given by
$$k = \sum_{ \alpha \in A_N, \ n\in \bbZ} b_{\alpha, n/2} c_{\alpha, -n}$$
where $b_{\alpha, n}$ are the Fourier coefficients of the Eisenstein series with the constant term
$b_{0, 0} = -1/2$ and $b_{\alpha, 0} = 0$ for $\alpha \not= 0$.}

\medskip

In the following we study the divisors $\sum_{\alpha \in A_N, n<0} c_{\alpha , n} \calH(N)_{\alpha, n}$
where $c_{\alpha, n}$ depends only on the type of $\alpha$.
Recall that there are 1, 15, 15, 1, 20 and 12 elements in $A_N$ of types
$(00)$, $(0)$, $(1)$, $(10)$, $(3/2)$ and $(1/2)$, respectively (see Lemma \ref{types}). 
We consider vector valued modular forms 
\begin{equation}\label{6dim}
f_{00},\ f_{0},\ f_{1},\ f_{10},\ f_{3/2},\ f_{1/2}
\end{equation}
where each $f_t$ is the sum of the $f_{\alpha}$ as $\alpha$ runs through the elements of $A_N$ of type $t$.  
Then we obtain a 6-dimensional representation
$$\rho^*: \SL(2,\bbZ) \to \GL(V)$$
induced by (\ref{dualweil2}) where $V$ is a 6-dimensional vector space (for simplicity, we use the same notation $\rho^*$).  Then
a direct calculation shows that $\rho^*$ is given by

\begin{equation}\label{}
\rho^*(T) =
\begin{pmatrix}1&0&0&0&0&0
\\0&1&0&0&0&0
\\0&0&-1&0&0&0
\\0&0&0&-1&0&0
\\0&0&0&0&\sqrt{-1}&0
\\0&0&0&0&0&-\sqrt{-1}
\end{pmatrix}.
\end{equation}

\begin{equation}\label{}
\rho^*(S) = {-\sqrt{-1}\over 8}
\begin{pmatrix}1&1&1&1&1&1
\\15&-1&-1&15&3&-5
\\15&-1&-1&15&-3&5
\\1&1&1&1&-1&-1
\\20&4&-4&-20&0&0
\\12&-4&4&-12&0&0
\end{pmatrix}.
\end{equation}

\subsection{Lemma}\label{cusp}
{\it The dimension of the space of modular forms of weight $3$ and of 
type $\rho^{*}$ is $2$.  The dimension of
the space of Eisenstein forms of weight $3$ and of type $\rho^{*}$ is also $2$.  In particular
there are no non-zero cusp forms in the obstraction space. 
}

\begin{proof}
The dimension of the space of modular forms of weight $3$ and of 
type $\rho^{*}$ is given by
$$d +dk/12 -\alpha(e^{\pi \sqrt{-1} k/2} \rho^*(S)) -\alpha ((e^{\pi\sqrt{-1}k/3} \rho^*(ST))^{-1}) - \alpha(\rho^*(T))$$
(\cite{B2}, section 4, \cite{F}, Proposition 2.1).
Here $k=3$ is the weight, 
$$d=\dim \{ x \in V : \rho^*(-E)x = (-1)^kx\}$$ 
and 
$$\alpha (A)=\sum_{\lambda} \alpha$$ 
where $\lambda$ runs through all eigenvalues of $A$ and $\lambda = e^{2\pi\sqrt{-1}\alpha}$, $0\leq \alpha <1.$
Since $\rho^*(S^2) = \rho^*(-E)$ and $\rho^*(S^2)(e_{\alpha}) = -e_{\alpha}$ for any
$\alpha$, we have $d = 6$.
A direct calculation shows that $\alpha(e^{\pi \sqrt{-1} k/2} \rho^*(S)) = 3/2$, $\alpha ((e^{\pi\sqrt{-1}k/3} \rho^*(ST))^{-1}) = 2$ and
$\alpha(\rho^*(T)) = 2$.  Thus we have proved the first assertion.  It follows from
\cite{F}, Remark 2.2 that
the space of Eisenstein forms of weight 3 and of type $\rho^{*}$ is isomorphic to
$$\{ x \in V \ : \ \rho^*(T)(x) = x, \ \rho^*(-E) = (-1)^k x\}$$
which has dimension 2.  Therefore the second and hence the third assertions hold.
\end{proof}
\smallskip

We need Fourier coefficients of Eisenstein series of weight $3$ with respect to $\rho^*$.
Since $\rho^*$ is trivial on the principal congurence subgroup of level $4$, these Eisenstein
series are linear combiantions of the standard ones.  
For $(a_1, a_2) \in (\bbZ/N\bbZ)^2$, 
let $G_k^{(a_1,a_2)}(\tau, N)$ be the Eisenstein series of weight $k$ and level $N$ 
corresponding to $(a_1,a_2)$ (e.g. see \cite{F}). 
Then 
\begin{equation}\label{eisen}
(c\tau + d)^{-k} G_k^{(a_1,a_2)}((a\tau + b)/(c\tau+d), N) = G_k^{(a_1, a_2)A}(\tau, N)
\end{equation}
where $A =
\begin{pmatrix}a&b
\\c&d
\end{pmatrix}
\in \SL(2,\bbZ)$.
Note that 
$$G_k^{(-a_1,-a_2)}(\tau, N) = (-1)^k G_k^{(a_1,a_2)}(\tau, N).$$
Put
$E_{1} = G_{3}^{(0,1)}(\tau, 4)$, $E_{2} = G_{3}^{(1,0)}(\tau,  4)$,
$E_{3} = G_{3}^{(1,1)}(\tau,4)$, $E_{4} = G_{3}^{(1,2)}(\tau,4)$, 
$E_{5} = G_{3}^{(1,3)}(\tau,4)$, $E_{6} = G_{3}^{(2,1)}(\tau,4)$.
It follows from the equation (\ref{eisen}) that the actions of $T$ and $S$ on these forms are as follows:
$T$ fixes $E_1$ and sends
$$E_2 \to E_3 \to E_4 \to E_5 \to E_2, \quad E_6 \to -E_6,$$
and $S$ sends, up to $\tau^3$, 
$$E_1\to E_2 \to -E_1, \quad E_3 \to E_5 \to -E_3, \quad E_4 \to -E_6 \to -E_4.$$
An elementary calculation shows that
Eisenstein forms of weight 3 and of type $\rho^{*}$ are given by
$$f_{00} = a E_{1} + {i(a+b)\over 8}(E_{2} + E_{3} + E_{4}+E_5 ),$$
$$f_{0} = b E_{1} + {i(15a-b)\over 8}(E_{2} + E_{3} + E_{4} +E_5 ),$$
$$
f_{1} = {i(15a -b)\over 8}(E_{2} - E_{3} + E_{4} -E_5) +b E_6 ,$$
$$f_{10} = {i(a + b)\over 8}(E_{2} -E_{3} + E_{4} - E_5)+aE_6 ,$$
$$
f_{3/2} = {i(20a +4b)\over 8}(E_{2} - iE_{3} - E_{4} +iE_5),$$
$$f_{1/2} = {i(12a - 4b)\over 8}(E_{2} +iE_{3} -E_{4} -iE_5),$$
where $a, b$ are parameters.  On the other hand, the Fourier series of $G_k^{(a_1,a_2)}(\tau, N)$ are known (e.g. Freitag \cite{F}, \S 1) .   It follows that
the Fourier series of $E_i$ are given by :
$$E_{1} =  {i(2\pi)^{3} \over 2^8}\{ -i + 4iq + \cdot \cdot \cdot \},$$

$$E_{2} =  {i(2\pi)^{3} \over 2^7}\{q^{1/4} + 4q^{1/2} + 8q^{3/4} + 16q + \cdot \cdot \cdot \},$$

$$E_{3} =  {i(2\pi)^{3} \over 2^7}\{i q^{1/4} - 4q^{1/2} -8i q^{3/4} + 16q + \cdot \cdot \cdot \},$$
$$E_{4} =  {i(2\pi)^{3} \over 2^7}\{-q^{1/4} + 4q^{1/2} -8q^{3/4} + 16q + \cdot \cdot \cdot \},$$
$$E_5= {i(2\pi)^{3} \over 2^7}\{-i q^{1/4} - 4q^{1/2} +8i q^{3/4} + 16q + \cdot \cdot \cdot \},$$
$$E_6={i(2\pi)^{3} \over 2^7}\{2i q^{1/2} + 0\cdot q + \cdot \cdot \cdot \}.$$
Put $a = -{2^{7} \over (2\pi)^{3}}$ and $b=0$.  Then we have
$$f_{00} = -1/2 + 10q + \cdot \cdot \cdot \ ,\quad
f_{0} = 120 q + \cdot \cdot \cdot \ , \quad 
f_{1} =  30q^{1/2} + \cdot \cdot \cdot \ ,$$
$$f_{10} =  4q^{1/2} + \cdot \cdot \cdot \ ,\quad
f_{3/2} =  10q^{1/4} + \cdot \cdot \cdot \ ,\quad
f_{1/2} =  48q^{3/4} + \cdot \cdot \cdot \ .$$
Combining Lemma \ref{cusp} and Theorem \ref{freitag}, we have

\subsection{Theorem}\label{borcherds}
{\it There exists a holomorphic automorphic form $\Phi_k$ on $\calD(N)$ with respect to a subgroup of $\O(N)$ of finite index whose weight $k$ 
is $4$, $10$, $30$ or $48$, and whose zero divisor is the Heegner divisor 
$\calH(N)_{\kappa_N}$, 
$\calH(N)_{3/2}$, $\calH(N)_{1}$ or $\calH(N)_{1/2}$ respectively}.

\smallskip

\section{Additive liftings}\label{GBL}

For an even lattice $L$ of signature 
$(2,n)$, the {\it additive lifting} is a correspondence  from vector-valued modular forms of weight $k$ and type $\rho$ to
automorphic forms of weight $k + n/2 -1$ on $\calD(L)$ with respect to $\widetilde{\O}(L)$.
By using additive liftings, we construct a 5-dimensional space of automorphic forms on 
$\calD(N)$ with respect to $\Gamma_N = \widetilde{\O}(N)$.

Let $W$ be the subspace in $\bbC [A_N]$ defined in (\ref{5-dim}).  
Recall that $\O(q_N)$ acts
on $\bbC[A_N]$ which commutes with the action of $\SL(2,\bbZ)$.  Hence $\O(q_N)$
acts on $W$.  First we consider the following special vectors in $W$.
Let $\eta(\tau)$ be the Dedekind eta function.  Then
$$\eta (\tau + 1)^{18}= -\sqrt{-1} \cdot \eta(\tau)^{18},\quad 
\eta(-1/\tau)^{18} = -\sqrt{-1} \tau^{9} \cdot \eta(\tau)^{18}.$$
Therefore, for $\theta \in W$, $\eta(\tau)^{18} \cdot \theta$ is a vector-valued modular form of weight
9 and of type $\rho$.  
By applying additive lifting (\cite{B1}, Theorem 14.3), we have
an automorphic form $F_{\theta}$ on $\calD(N)$ of weight $10$ associated to $\eta(\tau)^{18}\theta$.
Let $A_k(\calD(N), \Gamma_N)$ be the space of automorphic forms of
$\calD(N)$ of weight $k$ with respect to $\Gamma_N$.
Recall that $\O(q_N)$ acts on $W$, and also acts naturally on $\calD(N)$.
Additive lifting is an $\O(q_N)$-equivariant map
\begin{equation}\label{GBlift}
W \to  A_{10}(\calD(N), \Gamma_N)
\end{equation}
by sending $\theta$ to $F_{\theta}$.
We denote by $\widetilde{W}$ the image of the map (\ref{GBlift}).

\subsection{Lemma}\label{5dim}
{\it $\widetilde{W}$ is a $5$-dimensional space on which $\O(q_N)$ acts irreducibly.}
\begin{proof}
It suffices to see that the map (\ref{GBlift}) is non-zero.  Then the assertion follows from the Schur's lemma. 
We use Theorem 14.3 in \cite{B1}.
We fix an orthogonal decomposition
$N = U(2)\oplus U(2) \oplus A_1\oplus A_1$.  Let $e_1, f_1$ be a basis of the first factor 
$U(2)$ with $e_1^2 =f_1^2 =0, \la e_1, f_1\ra = 2$, and let
$e_2, f_2$ be a basis of the second $U(2)$ with $e_2^2=f_2^2=0, \la e_2, f_2\ra =2$.  
Let $a_1$ be a basis of the third factor $A_1$ and $a_2$ a basis of the fourth factor.
Put $z = e_1$ and $z' = f_1/2$ and let $K = z^{\perp}/\bbZ z = U(2)\oplus A_1\oplus A_1
\subset N$.   Obviously $A_N = U(2)^*/U(2)\oplus K^*/K$.
For simplicity, we denote by $\overline{x}$ for $x \ \mod \ N \in A_N$.
We consider the Fourier expansion around $z$.
Since 
$$\eta(\tau)^{18} = q^{3/4} + \cdots ,$$
for $\theta = (c_{\alpha})\in W$,
the initial term of 
$$\eta(\tau)^{18} \theta = \sum_{\alpha \in A_N} e_{\alpha} \sum_{n \in \bbQ} c_{\alpha}(n) e^{2\pi\sqrt{-1}n\tau}$$ 
is 
$$\sum_{\alpha \in A_N,\ q_N(\alpha) =3/2} e_{\alpha} c_{\alpha} q^{3/4}.$$  
Now we consider a special vector $\lambda = e_2/2+f_2 +a_1/2 $.  Since 
$\la \lambda, \lambda \ra = 3/2 > 0$, $\lambda$ has positive inner products with all elements in the interior of
the Weyl chamber.  Also note that $\overline{\lambda}$ is of type $(3/2)$.  
We choose $\theta = (c_{\alpha}) \in W$ satisfying $c_{\overline{\lambda}} \not=0$
and $c_{\overline{e_1/2 + \lambda}}=0$ (for example, we take the isotoropic subspace
$I$ generated by $\overline{f_1/2}$ and $\overline{e_2/2}$, and consider $V$ generated by $I$ and $\kappa_N=\overline{a_1/2 + a_2/2}$.  Then the support of 
$\theta_V$ is
$$\{ \ \overline{(f_1+a_i)/2},\quad \overline{(e_2+a_i)/2},\quad \overline{(f_1+e_2+a_i)/2} \ \}_{i =1,2},$$
and hence $\theta_V$ satisfies the condition).
Now it follows from \cite{B1}, Theorem 14.3 that the Fourier coefficient of $e^{2\pi \sqrt{-1} \la \lambda, Z\ra }$ in 
the lifting $F_{\theta}$ of $\eta(\tau)^{18} \theta$ is equal to 
$$c_{\overline{\lambda}}(\lambda^2/2)\cdot e^{2\pi \sqrt{-1} \la \lambda, z'\ra} + 
c_{ \overline{e_1/2 + \lambda}}((e_1/2 + \lambda)^2/2)\cdot e^{2\pi\sqrt{-1}\la e_1/2 + \lambda, z'\ra} = 
c_{\overline{\lambda}}(3/4) =  c_{\overline{\lambda}} \not= 0.$$
Hence the lifting of $\eta(\tau)^{18} \theta$ is non-zero.
\end{proof}

\subsection{Theorem}\label{GB}
{\it Let $\theta_V \in W$ be as in Lemma $\ref{key}$.
Let $F_V$ be the additive lifting of $\eta(\tau)^{18} \cdot \theta_V$.  Then $F_V$ is an automorphic form on 
$\calD(N)$ of weight $10$ with respect to $\Gamma_N$.  Moreover $F_V$ vanishes exactly along 
$$\sum_{\alpha\in V, \ q_N(\alpha)=1} {\calH(N)}_{\alpha}$$ 
with multiplicity one.
}
\begin{proof}
Let $\alpha \in V$ with $q_N(\alpha) = 1$.
Recall that the reflection $t_{\alpha}$ is represented by a reflection $s_r$ acting on 
$\calD(N)$ where $r \in N$ with $r^2 = -4$ and $r/2 \ \mod \ N = \alpha$
(see the equations (\ref{ref2}), (\ref{ref})).
It now follows from Lemma \ref{key} and the $\O(q_N)$-equivariance of the additive lifting (\ref{GBlift}) that $F_V$ vanishes along 
${\calH(N)}_{\alpha}$ where $\alpha \in V$ with $q_N(\alpha)=1$.  
Therefore the product of fifteen $F_V$ has weight $150$ and vanishes along Heegner divisors $\calH(N)_{\alpha}$ ($\alpha \in V, q_N(\alpha)=1, 
\alpha \not=\kappa_N$) with at least
multiplicity 3 and along Heegner divisor $\calH(N)_{\kappa_N}$ with at least multiplicity $15$.

On the other hand, consider the automorphic forms $\Phi_4$, $\Phi_{30}$ of weight 4, 30 in Theorem \ref{borcherds}.
Then $\Phi_{4}^{15} \cdot \Phi_{30}^3$ has weight $150$ and vanishes 
along Heegner divisors $\calH(N)_{\alpha}$ ($\alpha \in V, q_N(\alpha)=1, 
\alpha \not=\kappa_N$) with exactly multiplicity 3 and along the Heegner divisor $\calH(N)_{\kappa_N}$ with exactly multiplicity $15$.
Then the ratio 
$$\prod_V F_V / (\Phi_{4}^{15} \cdot \Phi_{30}^3)$$ 
is a holomorphic automorphic form of weight $0$, and hence it is constant by Koecher principle.
\end{proof}

\subsection{Remark}\label{1dim}
Let $W_0$ be the 1-dimensional subspace in $\bbC [A_N]$ on which the character of $\SL(2,\bbZ)$ is given by $\chi_3$ (Lemma \ref{char2}).   Let $\theta_0$ be a generator of $W_0$.
Then by definition,
$$\rho(S)(\theta_0) = \sqrt{-1} \theta_0, \quad \rho(T)(\theta_0) = \sqrt{-1} \theta_0.$$
Hence $\eta(\tau)^{6} \cdot \theta_0$ is a vector-valued modular form of weight
3 and of type $\rho$.  
Moreover we see that
$t_{\kappa_V}(\theta_0) = -\theta_0$ where $t_{\kappa_V}$ is the reflection associated 
with $\kappa_V$.  Similar argument as above shows that the additive lifting $F_0$ 
of $\eta(\tau)^{6} \cdot \theta_0$ is an automorphic form of weight 4 and vanishes along $\calH(N)_{\kappa_N}$.  Moreover
by considering the ratio $F_0/\Phi_4$, 
the {\it additive lifting} $F_0$ coincides with the {\it Borcherds product} $\Phi_4$ in Theorem \ref{borcherds}.

\section{Automorphic forms on $\calD(M)$}\label{G}

In this section we show the existence of some automorphic forms on $\calD(M)$ 
by restricting the automorphic forms on $\calD(N)$ obtained in the previous sections.
First we fix an embedding of $M$ into $N$ as follows.   Fix a decomposition
$$N = U(2) \oplus U(2) \oplus A_1 \oplus A_1.$$
Let $a_1$ (resp. $a_2$) be a basis of the first (resp. second) component $A_1$ in the above decomposition.  Then we consider $M$ as a sublattice of $N$ generated by
$U(2) \oplus U(2)$ and $a_1-a_2$.  Note that $M$ is the orthogonal complement of 
$a_1 + a_2$ and is primitive in $N$.
This embedding induces an embedding $\calD(M) \subset \calD(N)$.

\subsection{Lemma}\label{embed}
{\it Any $\gamma \in \widetilde{\O}(M)$ can be extended to an isometry of $N$ acting trivially on $A_N$.  In other words,
$\widetilde{\O}(M) \subset \widetilde{\O}(N)$.}

\begin{proof}
This follows from Proposition 1.5.1 in Nikulin \cite{N1}.
\end{proof}

\subsection{Lemma}\label{restrict}
{\it The restriction of $\calH(N)_{1}$, $\calH(N)_{3/2}$ or 
$\calH(N)_{1/2}$ to $\calD(M)$ is $\calH(M)_{1}$,  
$2\calH(M)_{7/4}$ or $2\calH(M)_{3/4}$ respectively.}
\begin{proof}
For $r \in N$, 
write $r = r_1 + {m\over 2}(a_1+a_2)$ where $r_1 \in M^*$ and $m \in \bbZ$.
If $r/2 \in N^*$, then $r_1/2 \in M^*$.
Now assume $r^2= -4, r/2 \in N^*$ and $r \ \mod\ N \not= \kappa_N$.  
Then if $(r_1)^2\geq 0$, then the hyperplane
$r_1^{\perp}$ does not meet $\calD(M)$.  
If $(r_1)^2 < 0$, then $(r_1)^2 = -4 + m^2 = -4, -3$.
The case $m=1$ ($(r_1/2)^2 = -3/4$) does not occur because 
the values of $q_M$ are $0, 1, 3/4, 7/4 \ \mod \ 2\bbZ$.  Hence we have $(r_1)^2 = -4$ and $m=0$.  Thus the restriction of $\calH(N)_{1}$ to $\calD(M)$
is $\calH(M)_{1}$.  Similarly in case that $r/2 \in N^*$ and 
$r^2 = -2$ or $r^2 = -6$, then $m\not=0$ because the norm of any vector in $M$ is divided by $4$.  If $r^2=-6$ and $m=2$, then $r_1 \in M$ with $r_1^2 = -2$.  This is a contradiction.  Hence $(r_1)^2 = -1$ if $r^2 =-2$ and $(r_1)^2= -5$ if $r^2=-6$.
Note that the hyperplanes $a_1^{\perp}$ and $a_2^{\perp}$ on
$\calD(N)$ cut the same hyperplane $(a_1-a_2)^{\perp}$ on $\calD(M)$.  Hence the restriction of 
$\calH(N)_{3/2}$ or $\calH(N)_{1/2}$ to 
$\calD(M)$ is $2\calH(M)_{7/4}$ or $2\calH(M)_{3/4}$
respectively.
\end{proof}

Recall that the 2-elementary subgroup $A^{(2)}_M $ of $A_M$ together with $q_M$ is a $5$-dimensional quadratic space over $\bbF_2$ with the radical $\kappa_M$.

\subsection{Lemma}\label{heegner1}
(1) {\it Let $I$ be a $2$-dimensional isotoropic subspace of $A_N$ and let $V$ be the subspace
generated by $I$ and $\kappa_N$.  Then
the restriction of the Heegner divisor 
$$\sum_{\alpha\in V, \ q_N(\alpha)=1, \ \alpha\not= \kappa_N} {\calH(N)}_{\alpha}$$
to $\calD(M)$ is of the form
\begin{equation}\label{heegner3}
\sum_{\beta\in V_1, \ q_M(\beta)=1, \ \beta\not= \kappa_M} {\calH(M)}_{\beta}
\end{equation}
where $V_1$ is a $3$-dimensional subspace of $A^{(2)}_M $ such that $b_M|V_1 \equiv 0$.}

(2) {\it The Heegner divisor $(\ref{heegner3})$
contains exactly $7$ lines through three $0$-dimensional components in a line.}

\begin{proof}
(1) The same argument in the proof of Lemma \ref{restrict} shows that the projections of three non-isotoropic vectors in $V$ not equal to $\kappa_N$ generate a 3-dimensional subspace $V_1 \subset A^{(2)}_M $ satisfying
$b_M|V_1\equiv 0$. 

(2) Let $\beta_1, \beta_2, \beta_3$ be non-isotoropic vectors in $V_1$ not equal to 
$\kappa_M$ and let $I$ be the maximal isotoropic subspace of $V_1$.  
Since $\dim V_1 =3$ and $b_M|V_1 \equiv 0$, $\kappa_M$ is contained in $V_1$, that is,
$\beta_1 + \beta_2 + \beta_3 =\kappa_M.$
For each $\beta_i \ (i=1,2,3)$, there are exactly 7 non-zero isotoropic vectors perpendicular to $\beta_i$, and hence
there exist exactly three 2-dimensional isotorpic subspaces
$I, I_i', I_i''$ perpendicular to $\beta_i$.  Note that $I, I_i', I_i''$ contains a non-zero isotoropic vector $\beta_i+\kappa_M$.  Denote by $\ell, \ell_i', \ell_i''$ the corresponding lines
respectively.  Then three lines meet at one point corresponding to $\beta_i+\kappa_M$ and the Heegner divisor $\calH(M)_{\beta_i}$ contains $\ell, \ell_i', \ell_i''$.
\end{proof}

\medskip
\noindent
 Now, by restricting the automorphic forms $\Phi_{30}, \Phi_{10}, \Phi_{48}$ in Theorem \ref{borcherds} to $\calD(M)$, we have the following:

\subsection{Corollary}\label{borchreds2}
{\it There exists automorphic forms $\Psi_{30}$, $\Psi_{10}$ or $\Psi_{48}$ on $\calD(M)$ of
weight $30$, $10$ or $48$ whose zero divisor on $\calD(M)$ is $\calH(M)_{1}$, $2\calH(M)_{7/4}$ or $2\calH(M)_{3/4}$
respectively}.
\begin{proof}
The restriction to $\calD(M)$ is not identically zero, and hence it is an automorphic form
with respect to a subgroup of $\O(M)$ of finite index by Lemma \ref{embed}.  Now 
the assertion follows from Lemmas \ref{restrict}.
\end{proof}

\subsection{Remark}\label{}
The restriction of $\Phi_4$ to $\calD(M)$ is identically zero.

\medskip
\noindent
Let $\alpha \in A_M$ be of type $(7/4)$ and let $r$ be in $M$ with $\la r, r \ra = -4$ and $r/4 \ \mod \ M =\alpha$.
Then the reflection $s_r$ defined by
$$s_r(x) = x - {2\la x, r\ra \over \la r, r \ra} r = x + \la x, r/2\ra r$$
is contained in $\O(M)$.  Since $r/2 \ \mod \ M = \kappa_M$, the action of $s_r$ on $A_M$ is equal to $-1_{A_M}$ (see the proof of Lemma \ref{types2}).
Therefore $s_r$ is contained in $\Gamma_M$.  Moreover the set of fixed points of $s_r$ is the hyperplane $r^{\perp}$.   Hence
the projection $\calD(M) \to \calD(M)/\Gamma_M$ is ramified along
the Heegner divisor of type $(7/4)$ with ramification degree two.  
Finally let $\Psi_{24}$ be the square root of $\Psi_{48}$.
Thus we have proved the following Theorem.

\subsection{Theorem}\label{Geer}
{\it The zero divisors of $\Psi_{30}$, $\Psi_{10}$ or $\Psi_{24}$ on 
$\calD(M)/\Gamma_M$ 
is $\calH(M)_{1}/\Gamma_M$, $\calH(M)_{7/4}/\Gamma_M$ or $\calH(M)_{3/4}/\Gamma_M$ respectively.}

\subsection{Remark}\label{Geer2}
The automorphic froms $\Psi_{30}, \Psi_{10}$ or $\Psi_{24}$ are known 
as Siegel modular forms (see \cite{vG}).  $\Psi_{30}, \Psi_{10}$ or $\Psi_{24}$ coincides with the Siegel modular form with the same Humbert surface as its zero divisor.  This follows from Koecher principle.  For example, 
$\Psi_{10}$ coincides with the product $\prod \theta_m^2(\tau)$ of the square of even theta constants.
The divisor $\calH(M)_{1}/\Gamma_M$, $\calH(M)_{7/4}/\Gamma_M$ or 
$\calH(M)_{3/4}/\Gamma_M$ on 
$\overline{\calD(M)/\Gamma_M}$ 
consists of fifteen, ten or six irreducible components, and coincides with the Humbert surface $H_4$, $H_1$ or $H_5$, respectively (see Remark \ref{Humbert}).

\medskip

Finally we consider the restriction of the 5-dimensional space $\widetilde{W}$ on
$\calD(N)$ to $\calD(M)$.  
Let $(x_1:\cdots : x_6)$ be a homogenous coordinate of $\bbP^5$.
Then Igusa quartic $\calI$ and Segre cubic $\calS$ are given by the equations (\ref{igusa}), (\ref{segre}), respectively.
We note that Igusa quartic is also called Castelnuovo-Richimond quartic (see Dolgachev \cite{D}, page 478).

The symmetric group $\mathfrak{S}_6$ naturally acts on 
$\calI$ and $\calS$ as automorphisms.  
It is classically known that the dual variety of $\calI$ is $\calS$ (\cite{Ba}, Chap.V).

Recall that the boundary of Satake's compactification $\overline{\mathfrak{H}_2/\Gamma_2(2)}$ 
consists of fifteen 1-dimensional and 15 0-dimensional components.  The fifteen 1-dimensional componets are 
the fifteen lines on $\calI$ defined by
\begin{equation}\label{lines}
(a:a:b:b:-a-b:-a-b)
\end{equation}
and its permutations (see \cite{vG}, Theorem 4.1).  Each line contains three 0-dimensional
components.  
For example, the line defined by (\ref{lines}) contains three 0-dimensional components 
$$(1:1:1:1:-2:-2), \ (1:1: -2 :-2  : 1:1), \ (-2:-2:1:1:1:1).$$
The singular locus of $\calI$ is nothing but the union of fifteen lines.  

On the other hand, the fifteen differences $x_i - x_j \  (i\not=j)$ are modular forms of weight 2 each of which 
defines an irreducible component of $H_4$ (\cite{vG}, \S 8).  
Each divisor defined by $x_i -x_j$ contains three lines meeting at one point.
For example, the divisor defined by $x_1-x_2$ contains three lines
$$(a:a:b:b:-a-b:-a-b), \ (a:a:b:-a-b:b:-a-b),\ (a:a:b:-a-b:-a-b:b)$$
which meet at $(2:2:-1:-1:-1:-1)$.  Moreover the divisor defined by the following modular form of weight $6$
\begin{equation}\label{cubics}
(x_i - x_j)(x_k-x_l)(x_m-x_n)\quad  (\{i,j,k,l,m,n\} = \{1,\ldots, 6\}).
\end{equation}
contains 7 lines containing one of the three 0-dimensional components on a line.  
For example,
$(x_1 - x_2)(x_3-x_4)(x_5-x_6)$ contains 7 lines passing a 0-dimensional component on
the line $(a:a:b:b:-a-b:-a-b)$.  Combining this with Lemma \ref{heegner1}, we have the following Lemma.

\subsection{Lemma}\label{V}
{\it There exists a bijective correspondence between
the fifteen Heegner divisors 
given in $(\ref{heegner3})$ and the fifteen divisors defined by }(\ref{cubics}).

\medskip
Let $\Phi_4$ be the automorphic form of weight $4$ with the Heegner divisor $\calH(N)_{\kappa_N}$ (Theorem \ref{borcherds}).  The ratio $F_V/\Phi_4$ is a
holomorphic automorphic form of weight $6$ whose zero divisor is 
$$\sum_{\alpha\in V, \ q_N(\alpha)=1, \alpha\not=\kappa_N} {\calH(N)}_{\alpha}$$
(see Theorem \ref{GB}).

\subsection{Theorem}\label{restr2}
{\it The linear system obtained by the restriction of fifteen $F_V/\Phi_4$ to $\calD(M)$ coincide with the one defined by
the modular forms given in } (\ref{cubics})
\begin{proof}
Note that the restriction of each $F_V/\Phi_4$ to $\calD(M)$ is an automorphic form of weight 6 with respect to a subgroup 
$\Gamma$ of $\O(M)$ of finite index (Theorem \ref{borcherds}, Theorem \ref{GB}, Lemma \ref{embed}).
Under the above identification given in Lemma \ref{V},
this form and the corresponding automorphic form $(\ref{cubics})$ have the same weight $6$ and the same zero divisor (Theorem \ref{borcherds}, Theorem \ref{GB}).  Hence they coinside by Koecher principle.  
\end{proof}
 
Finally we discuss the geometric meaning of the linear system $\calL$ generated by 15 cubics given in $(\ref{cubics})$.  Let $\bbP^4$ be the subspace of $\bbP^5$ defined by $\sum x_i =0$.
Consider $\calL$ as a linear system on $\bbP^4$.
Its base locus consists of the 15 lines defined by $x_i=x_j=x_k=x_l$.
For $(x_i) \in {\bbP}^5$, if we consider $x_i$ as an affine coordinate of $\bbP^1$, then fifteen functions
given in $(\ref{cubics})$ induce an $\mathfrak{S}_6$-equivariant isomorphism from the moduli space $P_1^6$ of ordered six points of $\bbP^1$ onto Segre cubic $\calS$ (see \cite{D}, Theorem 9.4.10, \cite{DO}, page 15).  This implies that the linear system $\calL$ on $\bbP^4$ defines an $\mathfrak{S}_6$-equivariant rational map $\Phi$ from $\bbP^4$ to $\calS$ whose general fiber is a rational curve.  The proof of the following theorem is due to I. Dolgachev.

\subsection{Theorem}\label{Igor}
{\it The linear system $\calL$  gives an $\mathfrak{S}_6$-equivariant rational map $\phi$  from $\calI$ to $\calS$ of degree $16$}.
\begin{proof}
The base locus of the rational map $\Phi : \bbP^4 \to \calS$ consists of 15 lines
through two points in the set of six points
\begin{equation}\label{}
\{ p_i = (x_1: \cdots : x_6) \ :  \ x_i = -5, \ x_j = 1\ (j\not=i) \}_{i=1,\ldots , 6}.
\end{equation}
Note that these six points are in general position.
Take a point $p \in \bbP^4$ such that 7 points $p, p_1,\ldots, p_6$ are in general position.  Then there is a unique rational normal curve $C$
of degree $4$ passing through $p, p_1,\ldots, p_6$ (\cite{GHa}, p.179,\ p.530).  
Each cubic in $\calL$ has singularities at six points $p_i$.  Therefore if $p$
is contained in a cubic in $\calL$, then $C$ is contained in this cubic.  This implies
that $C$ is a fiber of the map $\Phi$.  
Recall that the singular locus of $\calI$ is the union of fifteen lines (\cite{vG}, Theorem 4.1).
Since the six points $p_1,\ldots, p_6$ do not lie on Igusa quartic $\calI$, a general $C$ intersects $\calI$ at $4\times 4=16$ points.  Hence the resrtiction $\phi$ of $\Phi$ to $\calI$ has degree 16. 
\end{proof}

Since the dual variety of Segre cubic is Igusa quartic (\cite{Ba}), we have the following Corollary.

\subsection{Corollary}\label{selfmap}
{\it The rational map $\phi : \calI \to \calS$ induces a rational self-map of $\calI$ of
degree $16$.}

\subsection{Remark}\label{}
The author \cite{K3} gave a 5-dimensional linear system of holomorphic automorphic forms on a 3-dimensional complex ball by applying
Borcherds theory of automorphic forms.  This linear system gives the dual map from Segre cubic to Igusa quartic.  The author does not know the geometric meaning of the space 
$\widetilde{W}$ of automorphic forms on $\calD(N)$.


\end{document}